\documentclass{amsart}
\usepackage{amssymb,latexsym}

\newtheorem{thm}{Theorem}[section]
\newtheorem{prop}[thm]{Proposition}
\newtheorem{cor}[thm]{Corollary}
\newtheorem{lem}[thm]{Lemma}
\newtheorem{conj}[thm]{Conjecture}
\newtheorem{exa}[thm]{Example}

\newcommand{\ben}{\begin{enumerate}}
\newcommand{\een}{\end{enumerate}}
\newcommand{\ble}{\begin{lem}}
\newcommand{\ele}{\end{lem}}
\newcommand{\bth}{\begin{thm}}
\renewcommand{\eth}{\end{thm}}
\newcommand{\bpr}{\begin{prop}}
\newcommand{\epr}{\end{prop}}
\newcommand{\bco}{\begin{cor}}
\newcommand{\eco}{\end{cor}}
\newcommand{\bcon}{\begin{conj}}
\newcommand{\econ}{\end{conj}}
\newcommand{\bde}{\begin{defn}}
\newcommand{\ede}{\end{defn}}
\newcommand{\bex}{\begin{exa}}
\newcommand{\eex}{\end{exa}}
\newcommand{\barr}{\begin{array}}
\newcommand{\earr}{\end{array}}
\newcommand{\btab}{\begin{tabular}}
\newcommand{\etab}{\end{tabular}}
\newcommand{\beq}{\begin{equation}}
\newcommand{\eeq}{\end{equation}}
\newcommand{\bea}{\begin{eqnarray*}}
\newcommand{\eea}{\end{eqnarray*}}
\newcommand{\bce}{\begin{center}}
\newcommand{\ece}{\end{center}}
\newcommand{\bpi}{\begin{picture}}
\newcommand{\epi}{\end{picture}}
\newcommand{\bfi}{\begin{figure} \begin{center}}
\newcommand{\efi}{\end{center} \end{figure}}

\newcommand{\bsl}{\begin{slide}{}}
\newcommand{\esl}{\end{slide}}

\newcommand{\bib}{thebibliography}

\newcommand{\hs}[1]{\hspace{#1}}
\newcommand{\hso}[1]{\hspace{-1pt}}

\newcommand{\qmq}[1]{\quad\mbox{#1}\quad}

\newcommand{\rp}[2]{\rule{#1pt}{#2pt}}

\newcommand{\iso}{\cong}

\newcommand{\zh}{\hat{0}}
\newcommand{\oh}{\hat{1}}

\newcommand{\ptn}{\vdash}

\newcommand{\jn}{\vee}

\newcommand{\mt}{\wedge}

\newcommand{\case}[4]{\left\{\barr{ll}#1&\mbox{#2}\\#3&\mbox{#4}\earr\right.}

\def\<{\langle}
\def\>{\rangle}

\newcommand{\spn}[1]{\langle{#1}\rangle}
\newcommand{\ree}[1]{(\ref{#1})}

\newcommand{\ra}{\rightarrow}

\newcommand{\al}{\alpha}
\newcommand{\be}{\beta}

\newcommand{\de}{\delta}

\newcommand{\la}{\lambda}

\newcommand{\om}{\omega}

\newcommand{\rhot}{\tilde{\rho}}
\newcommand{\si}{\sigma}

\newcommand{\ze}{\zeta}

\newcommand{\La}{\Lambda}

\newcommand{\bx}{{\bf x}}

\newcommand{\bbQ}{{\mathbb Q}}

\newcommand{\bbZ}{{\mathbb Z}}

\newcommand{\cG}{{\mathcal G}}

\newcommand{\cM}{{\mathcal M}}

\newcommand{\cP}{{\mathcal P}}

\newcommand{\fS}{{\mathfrak S}}

\newcommand{\dil}{\displaystyle}

\newcommand{\aim}{Adv.\ in Math.\/}

\newcommand{\oup}{Oxford University Press}

\newcommand{\Ad}{\dot{A}}
\newcommand{\Add}{\ddot{A}}
\newcommand{\Addd}{A^{(n)}}
\newcommand{\xd}{\dot{x}}
\newcommand{\xdd}{\ddot{x}}
\newcommand{\xddd}{x^{(n)}}
\newcommand{\Td}{{\dot{T}}}
\newcommand{\Ud}{{\dot{U}}}
\newcommand{\bxd}{{\bf\dot{x}}}
\newcommand{\bxdd}{{\bf\ddot{x}}}
\newcommand{\bxddd}{{\bf x^{(n)}}}
\newcommand{\byd}{{\bf\dot{y}}}
\newcommand{\bydd}{{\bf\ddot{y}}}
\newcommand{\byddd}{{\bf y^{(n)}}}
\newcommand{\vla}{{\vec{\la}}}
\newcommand{\oned}{\dot{1}}
\newcommand{\onedd}{\ddot{1}}
\newcommand{\twod}{\dot{2}}
\newcommand{\twodd}{\ddot{2}}
\newcommand{\thrd}{\dot{3}}
\newcommand{\thrdd}{\ddot{3}}
\newcommand{\foud}{\dot{4}}
\newcommand{\foudd}{\ddot{4}}

\newcommand{\isdom}{\unlhd}
\newcommand{\vm}{\vec{m}}
\newcommand{\beh}{\hat{\be}}
\newcommand{\bec}{\check{\be}}

\begin{document}

\title{Symmetric Functions in Noncommuting Variables
}
\author{Mercedes H. Rosas}
\address{Departamento de Matem\'aticas,
Universidad Sim\'on Bol\'{\i}var, Apdo.\ Postal 89000,
Caracas, VENEZUELA}
\email{mrosas@ma.usb.ve}
\author{Bruce E. Sagan}
\address{Department of Mathematics, Michigan State University,
East Lansing, MI 48824-1027, USA}
\email{sagan@math.msu.edu}

\date{\today}
\keywords{noncommuting variables, partition lattice, Schur function,
  symmetric function}
\subjclass{Primary  05E05; Secondary 05E10, 05A18
}
	\begin{abstract}
Consider the algebra $\bbQ\spn{\spn{x_1,x_2,\ldots}}$ of formal power
series in countably many noncommuting variables over the rationals.
The subalgebra $\Pi(x_1,x_2,\ldots)$ of symmetric functions in
noncommuting variables consists 
of all elements invariant under permutation of the 
variables and of bounded degree.  We develop a theory of such
functions analogous to the ordinary theory of symmetric functions.  In
particular, we define analogs of the monomial, power sum, elementary,
complete homogeneous, and Schur symmetric functions as will as
investigating their properties.
	\end{abstract}
\maketitle

\section{Introduction}				
\label{i}

Let $\bbQ[[x_1,x_2,\ldots]]=\bbQ[[\bx]]$ be the algebra of formal
power series over $\bbQ$ in a countably infinite set of commuting
variables $x_i$.  For each positive integer $m$, the symmetric group
$\fS_m$ acts on $\bbQ[[\bx]]$ by 
\beq
\label{action}
g f(x_1,x_2,\ldots)=f(x_{g(1)},x_{g(2)},\ldots)
\eeq
where $g(i)=i$ for $i>m$.  We say that $f$ is {\it symmetric\/} if it is
invariant under the action of $\fS_m$ for all $m\ge1$.  The {\it
algebra of symmetric functions}, $\La=\La(\bx)$, consists of all
symmetric $f$ of bounded degree.  This algebra has a
long, venerable history  in 
combinatorics, algebraic geometry, and representation theory; see,
e.g., \cite{ful:yt, mac:sfh, sag:sg, sta:ec2}.

Now consider 
$\bbQ\spn{\spn{x_1,x_2,\ldots}}=\bbQ\spn{\spn{\bx}}$, the associative
algebra of formal power series in the {\it noncommuting\/} variables
$x_1,x_2,\ldots$  
Define the {\it algebra of symmetric functions
in noncommuting variables}, $\Pi=\Pi(\bx)$, to be the subalgebra
consisting of all elements in $\bbQ\spn{\spn{\bx}}$ which are invariant under
the action defined by~\ree{action} and 
of bounded degree.  (This is not to be confused with the algebra of
noncommutative symmetric functions of 
Gelfand et.\ al.~\cite{gkllrt:nsf} or 
the partially commutative symmetric functions studied by Lascoux and
Sch\"utzenberger~\cite{ls:mp} as well as by Fomin and Greene~\cite{fg:nsf}.)
This algebra was first studied by M. C. Wolf~\cite{wol:sfn} in 1936.
Her goal was to provide an analogue of the fundamental theorem of
symmetric functions in this context.  The concept then lay dormant for
over 30 years until Bergman and Cohn generalized Wolf's
result~\cite{bc:sef}.   Still later, Kharchenko~\cite{kha:aif} proved
that if $V$ is a graded vector space and $G$  a group of
grading-preserving automorphisms of the tensor algebra of $V$, then the
algebra of invariants of $G$ is also a tensor algebra.
Anick~\cite{ani:hit} then removed the condition
that $G$ preserve the grading.
Most recently, Gebhard
and Sagan~\cite{gs:csf} revived these ideas as a tool
for studying Stanley's chromatic symmetric function of a
graph~\cite{sta:sfg,sta:gcr}. 

This paper gives the first systematic study of $\Pi(\bx)$ and is
structured as follows.  In the next section we define
$\Pi(\bx)$-analogues for the monomial, power sum, elementary, and
complete homogeneous bases of $\La(\bx)$.  We relate the two sets of
bases in Theorem~\ref{rho} using the projection map
$\rho:\bbQ\spn{\spn{\bx}}\ra\bbQ[[\bx]]$ which lets the variables commute.
In section~\ref{cb}, we derive change of basis equations for these
four bases by summation or M\"obius inversion over
the lattice of partitions.  As a consequence, we obtain some
properties of the fundamental involution $\omega:\Pi(\bx)\ra\Pi(\bx)$ in
Theorem~\ref{om}.  In the following section, we define a right inverse for
$\rho$, called the lifting map, and study its relation with an inner product
on $\Pi(\bx)$.  In Section~\ref{msf}, we recall some facts about the
algebra of MacMahon
symmetric functions $\cM$ and show that a particular subspace of $\cM$
is naturally isomorphic to $\Pi(\bx)$ as a vector space.  This
permits us to define a noncommuting-variable analogue, $S_\la$, of a
Schur function in Section~\ref{sf}.  The next two sections are
devoted to obtaining analogues for $S_\la$ of the Jacobi-Trudi determinants
(Theorem~\ref{jt}) and Robinson-Schensted-Knuth algorithm 
(Theorem~\ref{rsk}).  We end with a list of comments and
open questions.  

\section{Basic definitions}
\label{bd}

Let $n$ be a nonnegative integer.
For $\la=(\la_1,\la_2,\ldots,\la_l)$ a
partition of $n$, we
write $\la\vdash n$ and denote the  length
of $\la$ by $l=l(\la)$.  We will  use the notation
\beq
\label{lam}
\la=(1^{m_1},2^{m_2},\ldots,n^{m_n})
\eeq
to mean that $i$ appears in $\la$ with multiplicity $m_i$, $1\le i\le n$.
The bases of the symmetric function algebra $\La(\bx)$ are indexed by
partitions.  Following~\cite{mac:sfh,sta:ec2},  we use the notation
$m_\la$, $p_\la$, $e_\la$, and $h_\la$ 
for the monomial, power sum, elementary, and complete
homogeneous symmetric functions bases.  Our next goal is to define the
analogues of these bases in noncommuting variables; these analogues
will be labeled by {\it set partitions}.

Define $[n]=\{1,2,\ldots,n\}$.
A {\it set partition $\pi$ of $[n]$} is a family
of disjoint sets, called {\it blocks} $B_1,B_2,\ldots, B_l$, whose
union is $[n]$.  We write $\pi=B_1/B_2/\ldots/B_l\vdash[n]$ 
and define {\it length\/} $l=l(\pi)$ as the number of blocks.
There is a natural mapping from set partitions to integer partitions
given by
$$
\la(\pi)=\la(B_1/B_2/\ldots/B_l)=(|B_1|,|B_2|,\ldots,|B_l|)
$$
where we assume that $|B_1|\ge|B_2|\ge\ldots\ge|B_l|$.
The integer partition $\la(\pi)$ is the {\it type\/} of the set
partition $\pi$.  

The partitions of $[n]$ form the {\it partition lattice\/}
$\Pi_n$.  (Do not confuse $\Pi_n$ with the algebra $\Pi(\bx)$.)
In $\Pi_n$ the ordering is by refinement: $\pi\le\si$ if
each block $B$ of $\pi$ is contained in some block $C$ of $\si$.  The 
meet (greatest lower bound) and join (least upper bound) operations in
$\Pi_n$ 
will be denoted $\mt$ and $\jn$, respectively.  There is a  rank
function in $\Pi_n$ given by $r(\pi)=n-l(\pi)$.

To obtain analogues of the bases of $\La(\bx)$ in this setting, 
it will be helpful to think of the elements
of $[n]$ as indexing the positions in a monomial 
$x_{i_1}x_{i_2}\cdots x_{i_n}$.  This makes sense because the
variables do not commute.  Now given $\pi\ptn[n[$, define the {\it
monomial symmetric function, $m_\pi$ in noncommuting variables\/} by
$$
\barr{ll}
m_{\pi}
= \dil\sum_{(i_{1},i_{2},\ldots ,i_{n})} x_{i_{1}}x_{i_{2}}\cdots x_{i_{n}}
&\mbox{where the sum is over all $n$-tuples   $(i_1,i_2,\ldots,i_n)$ with}
\\[-10pt]
&\mbox{$i_j=i_k$ if and only if $j,k$ are in the same
  block in $\pi$.} 
\earr
$$
For example, 
$$
m_{13/24}=x_1x_2x_1x_2 +x_2x_1x_2x_1  +x_1x_3x_1x_3+ 
x_3x_1x_3x_1+x_2x_3x_2x_3+x_3x_2x_3x_2  +\cdots
$$ 
These functions are precisely the symmetrizations of
monomials and so they are invariant under the action of $\fS_m$
defined previously. 
It follows easily that they form a basis for $\Pi(\bx)$.

We define the
{\it power sum function in noncommuting variables}, $p_\pi$,
by
$$
\mbox{
$p_{\pi}= 
\dil\sum_{(i_{1},i_{2},\ldots ,i_{n})} x_{i_{1}}x_{i_{2}}\cdots x_{i_{n}}$
where $i_j=i_k$ if $j,k$ are in the same block in $\pi$.
}
$$
To illustrate,
$$
p_{13/24}=x_1x_2x_1x_2 +x_2x_1x_2x_1+x_1^4+x_2^4+\cdots=m_{13/24}+m_{1234}.
$$

The {\it elementary symmetric function in noncommuting
variables} is
$$
\mbox{
$e_{\pi}= 
\dil\sum_{(i_{1},i_{2},\ldots ,i_{n})} x_{i_{1}}x_{i_{2}}\cdots x_{i_{n}}$
where $i_j\neq i_k$ if $j,k$ are in the same block in $\pi$.
}
$$
By way of example,
\bea
e_{13/24}&=&x_1x_1x_2x_2+x_2x_2x_1x_1+x_1x_2x_2x_1+x_2x_1x_2x_1+\cdots\\
	&=&m_{12/34}+m_{14/23}+m_{12/3/4}+m_{14/2/3}+m_{1/23/4}+m_{1/2/34}
		+m_{1/2/3/4}.
\eea

To define the analogue of the complete homogeneous symmetric functions, it will
be useful to introduce another way of looking at the previous
definitions.  (This was the method that Doubilet~\cite{dou:sft} used to
define certain ordinary symmetric functions associated with set
partitions.)  Any two sets $D,R$ and a function $f:D\ra R$ determine a {\it
kernel} set partition, $\ker f\vdash D$, whose blocks are the nonempty
preimages $f^{-1}(r)$ for $r\in R$.  For $f:[n]\ra \bx$ we denote by
$M_f$ the corresponding {\it monomial} 
$$
M_f=f(1)f(2)\cdots f(n).
$$
Directly from these definitions it follows that
$$
m_\pi=\sum_{\ker f=\pi} M_f.
$$
Using our running example, if $\pi=13/24$ then the functions with 
$\ker f=\pi$ are exactly those of the form $f(1)=f(3)=x_i$ and
$f(2)=f(4)=x_j$ where $i\neq j$.  This $f$ gives
rise to the monomial $M_f=x_ix_jx_ix_j$ in the sum for $m_{13/24}$.

Now define
\beq
\label{hpi}
h_\pi=\sum_{(f,L)} M_f
\eeq
where $f:[n]\ra \bx$ and $L$ is a linear ordering of the elements of
each block of $(\ker f)\mt\pi$.  Continuing with our running example,
$$
\barr{l}
h_{13/24}=m_{1/2/3/4}+m_{12/3/4}+2m_{13/2/4}+m_{14/2/3}
			+m_{1/23/4}+2m_{1/24/3}+m_{1/2/34}\\
\quad +m_{12/34}+4m_{13/24}+m_{14/23}
+2m_{123/4}+2m_{124/3}+2m_{134/2}+2m_{1/234}+4m_{1234}.
\earr
$$

Now we would like to give some justification to the above
nomenclature by exhibiting its relation to that used for the
ordinary symmetric  functions.  To this end, consider the
{\it projection\/} map
$$
\rho:\bbQ\spn{\spn{\bx}}\ra\bbQ[[\bx]]
$$
which merely lets the variables commute.  We will need the
notation
\bea
\la!	&=&\la_1!\la_2!\cdots\la_l!\\
\la^!	&=&m_1!m_2!\cdots m_n!
\eea
where the $m_i$ are the multiplicities in~\ree{lam}.  We extend
these conventions to set partitions by letting $\pi!=\la(\pi)!$ and
$\pi^!=\la(\pi)^!$.  Note that
\beq
\label{type}
{n\choose\la}:=\mbox{number of $\pi$ of type $\la$}=\frac{n!}{\la!\la^!}.
\eeq
The next proposition was
proved by Doubilet for his set partition analogues of ordinary
symmetric functions, and a similar proof can be given in the
noncommuting case.  The alternative demonstration given below
brings out the combinatorics behind some of Doubilet's
algebraic manipulations.
\bth
\label{rho}
The images of our bases under the projection map are:
\ben
\item[(i)] $\rho(m_\pi)=\pi^!  m_{\la(\pi)}$,
\item[(ii)] $\rho(p_\pi)= p_{\la(\pi)}$,
\item[(iii)] $\rho(e_\pi)=\pi! e_{\la(\pi)}$,
\item[(iv)] $\rho(h_\pi)=\pi! h_{\la(\pi)}$.
\een
\eth
\proof
For (i), let $B_1,B_2,\ldots,B_k$ be all blocks
of $\pi$ of a given size.  Recalling the remarks immediately before the
definition of $m_\pi$, we see that
$m_\pi$ is constant on the positions
indexed by each of these $B_i$.  Since $B_1, B_2,\ldots,B_k$  are of
the same size, the 
variables in the positions indexed by $B_i$ can be interchanged with
the variables in the positions indexed by $B_j$ for $1\le i,j\le k$ to give
another monomial in the sum for $m_\pi$ which maps to the same
monomial in the projection.  It follows that these blocks give
rise to a factor of $k!$ in the projection, and so $\pi$ will
contribute $\pi^!$.

\smallskip

To prove (ii), note that $p_{[n]}=m_{[n]}$ and so, from (i),
$\la(p_{[n]})=m_n=p_n$.  
Now $p_\la=p_{\la_1}p_{\la_2}\cdots p_{\la_l}$. Furthermore, 
$p_\pi=p_{B_1/B_2/\ldots/B_l}$ is precisely the shuffle of
$p_{[\la_1]}$, $p_{[\la_2]}$, $\ldots$, $p_{[\la_l]}$ where $|B_i|=\la_i$
and the elements from $p_{[\la_i]}$ are only permitted to be in the
positions indexed by $B_i$.  (We are letting the shuffle operation
distribute over addition.)  The desired equality follows.

\smallskip

The proof of (iii) is similar.  We have $\rho(e_{[n]})=n! e_n$ since if
all $n$ positions have different variables, then they can be permuted
in any of $n!$ ways and still give the same monomial in the
projection.  In the general case, we have the same phenomenon of
multiplication corresponding to shuffling, with each block $B$
contributing $|B|!$.  So the total contribution is $\pi!$.

\smallskip

Finally we consider (iv).  By the same argument as in (iii), it suffices
to show that $\rho(h_{[n]})=n! h_n$. Consider a monomial
$M=x_{j_1}^{\la_1}x_{j_2}^{\la_2}\cdots a_{j_l}^{\la_l}$ in $h_n$.
These variables can be rearranged to form $n!/\la!$
monomials in noncommuting variables where
$\la=(\la_1,\la_2,\ldots,\la_l)$.  To obtain one of these monomials 
in~\ree{hpi} we must have $\la(\ker f)=\la$
since $\ker f\mt [n]=\ker f$.  But then the number of pairs $(f,L)$ is
just $\la!$.  So the number of monomials in $h_\pi$ mapping to $M$
under $\rho$ is just $\la!\cdot n!/\la!=n!$, which is what we wanted.
\qed

We end this section by defining a second action of the symmetric group
which is
also interesting.  Since our variables do not commute, we can define
an action on places (rather than variables).  Explicitly, consider the
vector space of elements of $\Pi(\bx)$ which are homogeneous of
degree $n$.  Given a monomial of that degree, we define
\beq
\label{action2}
g\circ (x_{i_1}x_{i_2}\cdots x_{i_n})
=x_{i_{g(1)}}x_{i_{g(2)}}\cdots x_{i_{g(n)}} 
\eeq
and extend linearly.  It is easy to see that if $b_\pi$ is a basis
element for any of our four bases, then $g\circ b_\pi=b_{g\pi}$ where
$g$ acts on set partitions in the usual manner.

\section{Change of basis}
\label{cb}

We will now show that all the symmetric functions in noncommuting
variables defined in the previous section form bases for $\Pi(\bx)$.  Since
we already know this for the $m_\pi$, it suffices to find change of
basis formulas expressing each function in terms of the $m_\pi$
and vice-versa.  
Doubilet~\cite{dou:sft} has obtained
these results as well as those in the next section 
in a formal 
setting that includes ours as a special case.   But we replicate his
theorems and proofs here for completeness, to present them in standard
notation, and to extend and simplify some of them.

Expressing each symmetric function in terms of $m_\pi$ is easily done 
directly from the definitions, so the following proposition is given 
without proof.  In it, all lattice operations refer to $\Pi_n$ and
$\zh$ is the unique minimal element $1/2/\ldots/n$.
\bth
\label{btom}
We have the following change of basis formulae.
\ben
\item[(i)] $\displaystyle p_\pi=\sum_{\si\ge\pi} m_\si$,
\item[(ii)] $\displaystyle e_\pi=\sum_{\si\mt\pi=\zh} m_\si$,\rp{0}{20}
\item[(iii)] $\displaystyle h_\pi=\sum_\si (\si\mt\pi)! m_\si$.\rp{0}{20}\qed
\een 
\eth

To express $m_\pi$ in terms of the other functions, we will need the
M\"obius function of the parition lattice $\Pi_n$.  The {\it M\"obius
function\/} of any 
partially ordered set $P$ is the function $\mu:P\times P\ra\bbZ$
defined inductively by
$$
\mu(a,b)=\case{1}{if $a=b$,}{\displaystyle-\sum_{a\le c<b}\mu(a,c)}{else.}
$$
This can be rewritten in the useful and more intuitive form
\beq
\label{mu}
\sum_{a\le c\le b}\mu(a,c)=\de_{a,b}
\eeq
where $\de_{a,b}$ is the Kronecker delta.
For more information about M\"obius functions, see the seminal article
of Rota~\cite{rot:tmf} or the book of Stanley~\cite{sta:ec1}.  

The M\"obius function of $\Pi_n$ is well known.  In particular
$$
\mu(\zh,\oh)=(-1)^{n-1}(n-1)!
$$
where $\oh=12\ldots n$ is the unique maximal element of $\Pi_n$.
This is enough to determine $\mu$ on any interval of this lattice.
For example, for any $\pi=B_1/B_2/\ldots/B_l$ and $\la=\la(\pi)$ 
we have the lattice isomorphism $[\zh,\pi]\iso\prod_i \Pi_{\la_i}$.
Since the M\"obius function is preserved by isomorphism and
distributes over products, we have
$$
\mu(\zh,\pi)=\prod_i (-1)^{\la_i-1}(\la_i-1)!
$$
Note that, up to sign, this is just the number of permutations
$\al\in\fS_n$ which have disjoint cycle decomposition
$\al=\al_1\al_2\cdots\al_l$ where, for $1\le i\le l$, 
$\al_i$ is a cyclic permutation of the elements in $B_i$.  It follows that
$$
\sum_{\si\in\Pi_n} |\mu(\zh,\si)|=n!
$$
Or more generally, because of multiplicativity,
\beq
\label{sum||}
\sum_{\si\le\pi} |\mu(\zh,\si)|=\pi!
\eeq
a result which will be useful shortly.  We finally note that if
$\si=C_1/C_2/\ldots/C_m$ satisfies $\si\le\pi$ then we still  have
an isomorphism $[\si,\pi]\iso\prod_i \Pi_{\la_i(\si,\pi)}$ where
$\la(\si,\pi)$ is the integer partition whose $i$th part is the the
number of blocks of $\si$ contained in the $i$th block of $\pi$.
(We assume the blocks are listed so that the parts are in weakly
decreasing order).  Of course, $\la(\zh,\pi)$ is just the type of $\pi$.

The rest of the proofs in this section will all be based on the
M\"obius Inversion Theorem~\cite{rot:tmf,sta:ec1}.  We will also need
a simple corollary of 
that theorem which slightly generalizes a result of Doubilet~\cite{dou:sft}.
\bco
Let $P$ be a poset, let $F$ be a field, and consider three
functions  $f,g,h:P\ra F$ where $g(a)\neq 0$ for all $a\in P$. Then
$$\barr{l}
\displaystyle
f(a)=\sum_{b\le a} g(b)\sum_{c\ge b} h(c)\quad\mbox{for all $a\in P$}\\
\hs{80pt}\iff\displaystyle
h(a)=\sum_{c\ge a}\frac{\mu(a,c)}{g(c)}\sum_{b\le c}\mu(b,c)f(b)
\quad\mbox{for all $a\in P$}.
\earr$$
\eco
\proof
We will prove the forward direction as the converse is obtained by just
reversing the steps.  Doing (dual) M\"obius inversion on the outer
sum for $f(a)$ gives
$$
g(a)\sum_{c\ge a}h(c)=\sum_{b\le a}\mu(b,a)f(b).
$$
We can divide by $g(a)\neq0$ and then invert the sum
containing $h(c)$, which gives the desired result.\qed

We are now in a position to invert each of the sums in
Theorem~\ref{btom}.
\bth
We have the following change of basis formulae.
\ben
\item[(i$'$)] $\displaystyle m_\pi=\sum_{\si\ge\pi}\mu(\pi,\si) p_\si$,
\item[(ii$'$)] $\displaystyle m_\pi=
\sum_{\si\ge\pi}\frac{\mu(\pi,\si)}{\mu(\zh,\si)}
\sum_{\tau\le\si}\mu(\tau,\si) e_\tau$,\rp{0}{20}
\item[(iii$'$)] $\displaystyle m_\pi=
\sum_{\si\ge\pi}\frac{\mu(\pi,\si)}{|\mu(\zh,\si)|}
\sum_{\tau\le\si}\mu(\tau,\si) h_\tau$.\rp{0}{20}
\een 
\eth
\proof
Equation (i$'$) follows immediately from the M\"obius Inversion Theorem
applied to part (i) of Theorem~\ref{btom}.

\smallskip

For identity (ii$'$), use equation~\ree{mu} to write (ii) of
Theorem~\ref{btom} in the form
$$
e_\pi=\sum_\si\left(\sum_{\tau\le\si\mt\pi} \mu(\zh,\tau)\right) m_\si
     =\sum_{\tau\le\pi} \mu(\zh,\tau) \sum_{\si\ge\tau} m_\si.
$$
Using the corollary to invert this double sum gives the desired
result.

\smallskip

Finally consider (iii$'$).  Applying~\ree{sum||} to Theorem~\ref{btom} (iii)
gives
$$
h_\pi=\sum_\si\left(\sum_{\tau\le\si\mt\pi} |\mu(\zh,\tau)|\right) m_\si
     =\sum_{\tau\le\pi} |\mu(\zh,\tau)| \sum_{\si\ge\tau} m_\si.
$$
The corollary again provides the last step.\qed

The other bases-change equations are derived using similar techniques,
so we will content ourselves with merely stating the result after one
last bit of notation.
We define the {\it sign of
$\pi$}, $(-1)^\pi$, to be the sign of any permutation obtained by
replacing each block of $\pi$ by a cycle.  Note that
\beq
\label{sign}
\mu(\zh,\pi)=(-1)^\pi |\mu(\zh,\pi)|.
\eeq
\bth
\label{btob}
We have the following change of basis formulae.
$$
\barr{rclrcl}
e_\pi&=&\dil\sum_{\si\le\pi} \mu(\zh,\si)p_\si&\qquad
 p_\pi&=&\dil\frac{1}{\mu(\zh,\pi)}\sum_{\si\le\pi} \mu(\si,\pi)e_\si\\[25pt]
h_\pi&=&\dil\sum_{\si\le\pi} |\mu(\zh,\si)|p_\si&\qquad
 p_\pi&=&\dil\frac{1}{|\mu(\zh,\pi)|}\sum_{\si\le\pi} \mu(\si,\pi)h_\si\\[25pt]
e_\pi&=&\dil\sum_{\si\le\pi} (-1)^\si \la(\si,\pi)! h_\si&\qquad
	h_\pi&=&\dil\sum_{\si\le\pi} (-1)^\si \la(\si,\pi)! e_\si
\qed
\earr
$$
\eth

As an application of these equations, we will derive the properties of
an analogue of the involution $\omega:\La(\bx)\ra\La(\bx)$
defined by linearly extending $\omega(e_\la)=h_\la$.
Define a map on $\Pi(\bx)$, which we will also call $\omega$, by
$\omega(e_\pi)=h_\pi$
for all set partitions $\pi$ and linear extension.
\bth
\label{om}
The map $\omega:\Pi(\bx)\ra\Pi(\bx)$ has the following properties.
\ben
\item[(i)] It is an involution.
\item[(ii)] Each $p_\pi$ is an eigenvector for $\omega$ with
eigenvalue $(-1)^\pi$.
\item[(iii)] We have $\om\rho=\rho\om$.
\een
\eth
\proof
(i)  It suffices to show that the change of basis matrix between the
elementary and complete homogeneous symmetric functions equals its
inverse.  This follows directly from Theorem~\ref{btob}.

\smallskip

(ii) We merely compute the action of $\omega$ on the power sum basis
by expressing it in terms of the elementary symmetric functions and
using equation~\ree{sign}
$$
\omega(p_{\pi})
	=\omega\left(\frac{1}{\mu(\zh,\pi)} 
		\sum_{\sigma \le \pi} \mu(\sigma, \pi) e_{\sigma}\right)
	=\frac{1}{\mu(\zh,\pi)} 
		\sum_{\sigma \le \pi} \mu( \sigma, \pi) h_{\sigma}
	=\frac{(-1)^\pi}{|\mu(0,\pi)|} 
		\sum_{\sigma \le \pi} \mu( \sigma, \pi) h_{\sigma}
	=(-1)^\pi p_{\pi}.
$$

\smallskip

(iii)  It suffices to show that the desired equation holds on a basis.  So we
compute using Theorem~\ref{rho}~(iii) and~(iv) :
$\omega\rho(e_\pi)=\omega(\pi!e_{\la(\pi)})=\pi!h_{\la(\pi)}=\rho(h_\pi)
=\rho\omega(e_\pi)$.\qed

\section{The lifting map and inner products}
\label{lmip}

We will now introduce a right
inverse $\rhot$ for the projection map $\rho$ and an inner product for
which $\rhot$ is an isometry.  Define the {\it lifting map\/} 
$\rhot:\La(\bx)\ra\Pi(\bx)$
by linearly extending
\beq
\label{rhot}
\rhot(m_\la)=\frac{\la!}{n!}\sum_{\la(\pi)=\la} m_\pi.
\eeq
\bpr
The map $\rho\rhot$ is the identity map on $\La(\bx)$.
\epr
\proof Equation~\ree{type} and Theorem~\ref{rho}~(i) give
$\rho\rhot(m_\la)
	=\frac{\la!}{n!}\sum_{\la(\pi)=\la} \la^! m_\la=m_\la.$\qed

Recall that the standard inner product on $\La(\bx)$ is defined by
$\spn{m_\la,h_\mu}=\de_{\la,\mu}$.  We define its analogue in
$\Pi(\bx)$ by
\beq
\label{spn}
\spn{m_\pi,h_\si}=n!\de_{\pi,\si},
\eeq
where $\pi\vdash[n]$.  This bilinear form respects the
grading of $\Pi(\bx)$ in the sense that if $f,g$ are homogeneous
symmetric functions of different degrees, then $\spn{f,g}=0$.
\bth
The bilinear form $\spn{\cdot,\cdot}$  is symmetric, positive definite, and 
invariant under the action~\ree{action2}.
\eth

\proof
For symmetry, it suffices to show that
$\spn{h_\pi,h_\si}=\spn{h_\si,h_\pi}$.  By
Theorem~\ref{btom}~(iii),
\beq
\label{sym}
\spn{h_\pi,h_\si}
=\left\langle\sum_\tau (\pi\mt\tau)! m_\tau,h_\si\right\rangle
=n!(\pi\mt\si)!
\eeq
where we let $(\pi\mt\si)!=0$ if $\pi$ and $\si$ are partitions of
different sets.  Noting that  $(\pi\mt\si)!=(\si\mt\pi)!$ completes
the proof of symmetry. 

As for positive definiteness, take $f\in\Pi(\bx)$
and write $f=\sum_\pi c_\pi p_\pi$ for certain
coefficients $c_\pi$.  Then, using the expansions for the power sums
in Theorems~\ref{btom} and~\ref{btob}, we have 
\bea
\spn{f,f}
&=&\dil\left\langle\sum_\si c_\si p_\si,\ \sum_\tau c_\tau p_\tau
\right\rangle\\[10pt]
&=&\dil\left\langle\sum_\si c_\si \sum_{\pi\ge\si}m_\pi,\
\sum_\tau c_\tau \frac{1}{|\mu(\zh,\tau)|}\sum_{\pi\le\tau}\mu(\pi,\tau)h_\pi
\right\rangle\\[10pt]
&=&\dil n!\sum_\pi  \left(\sum_{\si\le\pi} c_\si\right)
\left(\sum_{\tau\ge\pi}c_\tau\frac{1}{|\mu(\zh,\tau)|}\mu(\pi,\tau)\right).
\eea
Now the coefficient of $c_\si c_\tau$ in this last sum is
$$
\frac{n!}{|\mu(\zh,\tau)|}\sum_{\si\le\pi\le\tau} \mu(\pi,\tau)
=\frac{n!\de_{\si,\tau}}{|\mu(\zh,\tau)|}.
$$
Since this is zero for $\si\neq\tau$ and positive otherwise, our form
is positive definite.

Finally, it suffices to verify invariance under the action on a pair of bases:
$$
\spn{g\circ m_\pi, g\circ h_\si} =\spn{m_{g\pi}, h_{g\si}} =
n!\de_{g\pi,g\si}=n!\de_{\pi,\si}=\spn{m_\pi,h_\si}.\qed
$$

\bth
The map $\rhot:\La(\bx)\ra\Pi(\bx)$ an isometry, i.e., 
$\spn{f,g}=\spn{\rhot(f),\rhot(g)}$ for $f,g\in\La(\bx)$. 
\eth
\proof It suffices to show that 
$\spn{\rhot(m_\la),\rhot(h_\mu)}=\spn{m_\la,h_\mu}$ for all $\la,\mu$.
To compute $\rhot(h_\mu)$, consider 
$$
H_\mu=\sum_{\la(\pi)=\mu} h_\pi.
$$
Expressing $H_\mu$ in terms of the monomial symmetric function basis
using Theorem~\ref{btom}~(iii), we see that the coefficient of $m_\si$
is the sum of $(\pi\mt\si)!$ over all $\la(\pi)=\mu$.  But the usual
action of the symmetric group on set partitions shows that this quantity
only depends on $\la(\si)$.  Thus by~\ree{rhot}, $H_\mu$ must be in
the image of $\rhot$.  Since $\rho$ is a left-inverse for $\rhot$, we
see that $H_\mu$ is the image under $\rhot$ of
$$
\rho(H_\mu)=\sum_{\la(\pi)=\mu} \rho(h_\pi)=
\sum_{\la(\pi)=\mu} \la(\pi)! h_{\la(\pi)}=
{n\choose\mu}\mu! h_\mu=\frac{n!}{\mu^!} h_\mu.
$$
So finally
$$
\spn{\rhot(m_\la),\rhot(h_\mu)}
=\left\langle\frac{\la!}{n!}\sum_{\la(\pi)=\la} m_\pi,
		\frac{\mu^!}{n!}\sum_{\la(\si)=\mu} h_\si\right\rangle
=\frac{\la!\la^!}{n!^2}{n\choose\la} n!\de_{\la,\mu}
=\de_{\la,\mu}
=\spn{m_\la,h_\mu}.
\qed
$$

To define the inner product~\ree{spn} in terms of other pairs of
bases, we will need the {\it zeta
function\/} of the partition lattice $\Pi_n$, defined by
$$
\zeta(\pi,\si)=\case{1}{if $\pi\le\si$,}{0}{else.}
$$
\bth
\label{bil}
The following formulae define equivalent bilinear forms.
$$
\barr{rclrcl}
\spn{e_\pi,e_\si}&=&n!(\pi\mt\si)!&\qquad
	\spn{e_\pi,h_\si}&=&n!\de_{\pi\mt\si,\zh}\\[10pt]
\spn{e_\pi,p_\si}&=&(-1)^\si n!\ze(\si,\pi)&\qquad
	\spn{e_\pi,m_\si}&=&(-1)^\si n!\la(\si,\pi)!\ze(\si,\pi)\\[10pt]
\spn{h_\pi,h_\si}&=&n!(\pi\mt\si)!&\qquad
	\spn{h_\pi,p_\si}&=&n!\ze(\si,\pi)\\[10pt]
\spn{h_\pi,m_\si}&=&n!\de_{\pi,\si}&\qquad
	\spn{p_\pi,p_\si}&=&n!\dil\frac{\de_{\pi,\si}}{|\mu(\zh,\pi)|}\\[10pt]
\spn{p_\pi,m_\si}&=&n!\dil\frac{\mu(\si,\pi)\ze(\si,\pi)}{|\mu(\zh,\pi)|}&
	\spn{m_\pi,m_\si}&=&
n!\dil\sum_{\tau\ge\pi\jn\si}\frac{\mu(\pi,\tau)\mu(\si,\tau)}{|\mu(\zh,\tau)|}
\qed
\earr
$$
\eth
The proof is similar to the derivation of~\ree{sym}, and is omitted.

\section{MacMahon symmetric functions}
\label{msf}

Schur functions in noncommuting variables will be defined in
Section~\ref{sf}.  This will require another piece of machinery,
namely the MacMahon symmetric functions.  The connection between
symmetric functions in noncommuting variables and MacMahon symmetric functions
was first pointed out by Rosas~\cite{ros:msf,ros:sms}.

Consider $n$ sets, each consisting of a countably infinite number of
commuting variables,
\bea
\bxd 	&=&\{\xd_1,\xd_2,\ldots\},\\
\bxdd	&=&\{\xdd_1,\xdd_2,\ldots\},\\
     	&\vdots&\\
\bxddd	&=&\{\xddd_1,\xddd_2,\ldots\}.
\eea
For each positive integer $m$, the symmetric group $\fS_m$ acts on
$\bbQ[[\bxd,\bxdd,\ldots,\bxddd]]$
diagonally, i.e.,
\beq
\label{mac}
g f(\xd_1,\xdd_1,\ldots,\xd_2,\xdd_2,\ldots)
= f(\xd_{g(1)},\xdd_{g(1)},\ldots,\xd_{g(2)},\xdd_{g(2)},\ldots).
\eeq
where $g(i)=i$ for $i>m$.  We say that
$f\in\bbQ[[\bxd,\bxdd,\ldots,\bxddd]]$ is {\it symmetric\/} if it is 
invariant under the action of $\fS_m$ for all $m\ge1$.

Consider a monomial
$$
M=\xd_1^{a_1}\xdd_1^{b_1}\cdots \left(x_1^{(n)}\right)^{c_1}
\xd_2^{a_2}\xdd_2^{b_2}\cdots \left(x_2^{(n)}\right)^{c_2}\cdots.
$$
Letting $\la^i=[a_i,b_i,\ldots,c_i]$ be the exponent sequence of the
variables of subscript $i$, we define the {\it multiexponent of $M$\/}
to be the vector partition
$$
\vla=\{\la^1,\la^2,\ldots\}
	=\{[a_1,b_1,\ldots,c_1],[a_2,b_2,\ldots,c_2],\ldots\}.
$$
By summing up the vectors which make up the parts of $\vla$ we
get the {\it multidegree of $M$\/}
$$
\vec{m}=[m_1,m_2,\ldots,m_n]
	=[a_1,b_1,\ldots,c_1]+[a_2,b_2,\ldots,c_2]+\cdots.
$$
In this situation we write $\vla\vdash\vec{m}$, $\vec{m}\vdash m$
where $m=\sum_i m_i$, and call $m$ the {\it degree of $M$}.  We say
that $f\in\bbQ[[\bxd,\bxdd,\ldots,\bxddd]]$ has {\it bounded
degree\/} if there is a positive integer $m$ such that all monomials
in $f$ have degree at most $m$.  Define
the {\it algebra of MacMahon symmetric functions},
$\cM=\cM(\bxd,\bxdd,\ldots,\bxddd)$, to be the subalgebra of
$\bbQ[[\bxd,\bxdd,\ldots,\bxddd]]$ consisting of all $f$ which are
symmetric under the action defined by~\ree{mac} and of bounded degree.

Given a vector partition $\vla$, there is an associated {\it monomial
MacMahon symmetric function\/} defined by
$$
m_\vla=\mbox{ sum of all the monomials with multiexponent $\vla$.}
$$
By way of example,
$$
m_{[2,1],[3,0]}=\xd_1^2\xdd_1\xd_2^3+\xd_1^3\xd_2^2\xdd_2+\cdots.
$$
Note that we drop the curly brackets around $\vla$ for readability.
These functions are precisely the symmetrizations of
monomials in $\bbQ[[\bxd,\bxdd,\ldots,\bxddd]]$ 
and so they are invariant under the action~\ree{mac} of $\fS_m$ for
all $m\ge1$.  It follows easily that they form a basis for $\cM$.

Call a basis $b_\vla$ of
$\cM$ {\it multiplicative\/} if it satisfies
$$
b_\vla=b_{\la^1}b_{\la^2}\cdots b_{\la^l}.
$$
We now define the bases of {\it power sum,
elementary, and complete homogeneous MacMahon symmetric functions\/}
to be multiplicative with
\bea
p_{[a,b,\ldots,c]} &=& m_{[a,b,\ldots,c]}\\[5pt]
\displaystyle
\sum_{a,b,\ldots,c} e_{[a,b,\ldots,c]} q^a r^b\cdots s^c
&=&
\displaystyle
\prod_{i\ge1} \left(1+\xd_i q + \xdd_i r +\cdots +\xddd_i s\right)\\[5pt]
\displaystyle
\sum_{a,b,\ldots,c} h_{[a,b,\ldots,c]} q^a r^b\cdots s^c
&=&
\displaystyle
\prod_{i\ge1} \frac{1}{1-\xd_i q - \xdd_i r -\cdots -\xddd_i s}.
\eea

To see the connection with noncommutative symmetric functions,
let $[1^n]$ denote the vector of $n$ ones.  Now consider
the subspace $\cM_{[1^n]}$ of $\cM$ spanned by all the $m_\vla$ where
$\vla\vdash[1^n]$.  There is a linear map 
$$
\Phi:\bigoplus_{n\ge0} \cM_{[1^n]}\ra\Pi
$$
given by
$$
\xd_i\xdd_j\cdots\xddd_k\stackrel{\Phi}{\mapsto}x_ix_j\cdots x_k.
$$
Given any $B\subseteq[n]$, the {\it characteristic vector of $B$\/} is
$[b_1,b_2,\ldots,b_n]$ where $b_i= 1$ if $i\in B$ and $b_i=0$ otherwise.
\bth[\cite{ros:msf}]
The map $\Phi$ is an isomorphism of vector spaces.  Furthermore, for
each basis we have discussed
$$
b_{\la^1,\la^2,\ldots,\la^l}\stackrel{\Phi}{\mapsto}b_{B_1/B_2/\ldots/B_l}
$$
where $b=m$, $p$, $e$, or $h$, and $\la^i$ is the characteristic vector
of $B_i$.\qed
\eth
By way of illustration
$b_{[1,0,1,0],[0,1,0,1]}\stackrel{\Phi}{\mapsto}b_{13/24}$ 
for any of our bases.

\section{Schur functions}
\label{sf}

We will now give a combinatorial definition of an analogue of a Schur
function in the 
setting of MacMahon symmetric functions.  This will give, via the map
$\Phi$, such a function in noncommuting variables. 
Consider the alphabet
\bea
A&=&\Ad\uplus\Add\uplus\cdots\uplus\Addd\\
&=&\{\oned,\twod,\ldots\}\uplus\{\onedd,\twodd,\ldots\}\uplus\cdots
	\uplus\{1^{(n)},2^{(n)},\ldots\}.
\eea
Partially order $A$ by 
\beq
\label{po}
i^{(k)}<j^{(l)} \qmq{if and only if} i<j.  
\eeq
Consider a partition $\la$ and a vector $\vec{m}=[m_1,\ldots,m_n]$
such that $\la,\vec{m}\ptn m$ for some nonnegative integer $m$.
Define a {\it dotted Young tableaux $\Td$ of  shape $\la$ and multidegree 
$\vec{m}=[m_1,\ldots,m_n]$\/} to be 
a filling of the shape of $\la$ (drawn in English style) with
elements of $A$ so that 
rows are nondecreasing, columns are strictly increasing, and there are
$m_k$ entries with $k$ dots.  
Now define the corresponding {\it MacMahon Schur function} to be
$$
S_\la^{\vec{m}}=\sum_{\la(\Td)=\la} M_\Td\qmq{where} 
M_\Td=\prod_{i^{(j)}\in\Td} x_i^{(j)}
$$
and the factor $x_i^{(j)}$ occurs in the above product with
multiplicity, i.e.,  the same
number of times that $i^{(j)}$ occurs in $\Td$.
For example, if $\la=(3,1)$ and $\vec{m}=[2,2]$, then the
coefficient of $\xd_1^2\xdd_1\xdd_2$ in $S_\la^{\vec{m}}$ is 3,
corresponding to the three dotted tableaux
$$
T_1=\barr{ccc} \oned&\oned&\onedd\\	\twodd	\earr,\quad
T_2=\barr{ccc} \oned&\onedd&\oned\\	\twodd	\earr,\quad
T_3=\barr{ccc} \onedd&\oned&\oned\\	\twodd	\earr.
$$
The notion of multidegree generalizes to any multiset $M$ of elements
from $A$.  If $M$ has $m_k$ elements with $k$ dots, we write
$\vm(M)=\vm=[m_1,\ldots,m_n]$.

\bth
The function $S_\la^{\vec{m}}$ is a MacMahon symmetric function.
\eth
\proof
It is obvious that $S_\la^{\vec{m}}$ is of bounded degree, so we need
only show that it is symmetric.
Because any permutation is a product of adjacent transpositions, it
suffices to show that $S_\la^{\vec{m}}$ is invariant under the
transposition $(i,i+1)$ where $i\ge1$.  So it suffices to find a
shape-preserving involution on  dotted tableaux
$\Td\ra\Td'$ which exchanges the number 
of elements equal to $i^{(k)}$ with the number equal to $(i+1)^{(k)}$
for all $k$, $1\le k\le n$.  We will use a generalization of a map of
Knuth~\cite{knu:pmg} used to prove that the ordinary Schur functions
are symmetric.

Since $\Td$ is semistandard, each column contains either a pair
$i^{(k)}, (i+1)^{(l)}$; exactly one of $i^{(k)}$ or $(i+1)^{(l)}$; or neither.
In the first case, replace the pair by $i^{(l)}, (i+1)^{(k)}$.  In the
second, replace $i^{(k)}$ by $(i+1)^{(k)}$ or replace $(i+1)^{(l)}$ by
$i^{(l)}$ as appropriate.  And in the third case there is nothing to
do.  It is easy to verify that this involution has the
desired properties. \qed

If $\vec{m}=[1^n]$ then we will write $S_\la$ for
$S^{\vec{m}}_\la$ and make no distinction between $S_\la$ and its image
under the map $\Phi$.  The latter will cause no problems because we
will never be multiplying these functions.  Note also that if
$\vec{m}$ has only one component, then $S^{\vec{m}}_\la=s_\la$, the
ordinary Schur function.

The $S_\la$ do not form a basis for $\Pi(\bx)$ since we
only have one for every integer, rather than set, partition.  
However, we can still provide analoques of some of the familiar
properties of ordinary Schur functions.
To state our  results, we will need the dominance order on integer
partitions, $\mu \isdom \la$, and the  Kostka numbers, $K_{\la,\mu}$.
(For definitions, see~\cite{mac:sfh,sta:ec2}.) 
\bth
\label{Sprop}
The functions $S_\la$ have the following properties.
\ben
\item[(i)] $\dil S_\la=
  \sum_{\mu\isdom\la} \mu!K_{\la,\mu}\sum_{\la(\si)=\mu} m_\si.$
\item[(ii)] The $S_\la$ are linearly independent.\rule{0pt}{15pt}
\item[(iii)] $\rho(S_\la)=n!s_\la$.
\item[(iv)] $\rhot(n!s_\la)=S_\la$. 
\item[(v)]  $\spn{S_\la,S_\mu}=n!^2\de_{\la,\mu}$.
\een
\eth
\proof
(i)  Consider a monomial $x^\Td$ where $\Td$ has shape $\la$ and
suppose that this monomial occurs in $m_\si$ where $\la(\si)=\mu$.
Then the number of ordinary tableaux $T$ with the same content, $\mu$,
as $\Td$ is $K_{\la,\mu}$ and this is only nonzero for $\mu\isdom\la$.
Furthermore, the number of ways to distribute dots in $T$ so as to
give the same monomial as $x^\Td$ is $\mu!$, so this finishes the
proof.

\smallskip

(ii) The lexicographic order on integer partitions is a linear
extension of the dominance order.
So from~(i), each $S_\la$ only contains $m_\si$ where $\la(\si)$ is
lexicographically less than or equal to $\la$, and those with
$\la(\si)=\la$ have nonzero coefficient.  So if one orders the $S_\la$
this way, then each Schur function will contain at least one 
monomial symmetric function not found previously in the list.

\smallskip

(iii) Using (i) again along with Theorem~\ref{rho}~(i) and
equation~\ree{type} gives 
$$
\rho(S_\la)=\sum_{\mu\isdom\la} \mu!K_{\la,\mu}\sum_{\la(\si)=\mu}\rho(m_\si)
=\sum_{\mu\isdom\la} \mu!K_{\la,\mu}\mu^!{n\choose \mu} m_\mu
=n!\sum_{\mu\isdom\la} K_{\la,\mu} m_\mu=n!s_\la.
$$

\smallskip

(iv) Clearly from~(i), all $m_\si$ with $\la(\si)=\mu$ have the same
coefficient in $S_\la$.  So $S_\la$ is in the image of $\rhot$.  The
equality now follows from~(iii) and the fact that $\rho$ is a
left-inverse for $\rhot$.

\smallskip

(v) We compute using~(iv) and the fact that $\rhot$ is an isometry
$$
\spn{S_\la,S_\mu}=\spn{\rhot(n!s_\la),\rhot(n!s_\mu)}=
\spn{n!s_\la,n!s_\mu}=n!^2\de_{\la,\mu}.\qed
$$

\section{Jacobi-Trudi determinants}
\label{jtd}

In this section, we prove analogs of the Jacobi-Trudi
determinants~\cite{mac:sfh} for the
$S_\la^{\vec{m}}$, where $\vec{m}$ is arbitrary.  The
ordinary and noncommuting variable cases are obtained as
specializations.  We use the lattice-path approach implicit in
Lindstr\"om~\cite{lin:vri} and developed explicitly by Gessel and
Viennot~\cite{gv:bdp}; see~\cite{sag:sg} for an exposition.

If $f\in \bbQ[[\bxd,\bxdd,\ldots,\bxddd]]$ and $\vec{m}$ is a vector,
then let $\spn{\vec{m}} f$ denote the sum of all terms $q_M M$ in $f$
where $q_M\in\bbQ$ and $M$ is a monomial of multidegree $\vec{m}$.
Also let $\la'$ denote the conjugate of the partition $\la$.

\bth 
\label{jt}
Given a partition $\la$ and vector $\vec{m}$ with
$\la, \vec{m}\vdash m$, we have
$$
S_\la^{\vec{m}}
=\spn{\vec{m}}\det\left(\sum_{\vec{t}\ \vdash\la_i-i+j} 
	h_{\vec{t}}\ \right)
$$
and
$$
S_{\la'}^{\vec{m}}
=\spn{\vec{m}}\det \left(\sum_{\vec{t}\ \vdash\la_i-i+j} 
	e_{\vec{t}}\ \right).
$$
\eth
\proof
We will only prove the first identity as the second is obtained by
a similar argument.

Consider infinite paths in the extended integer lattice
$\bbZ\times(\bbZ\uplus\infty)$:
$$
P= s_1, s_2, s_3,\ldots
$$
where the $s_t$ are steps of unit length either northward (N) or
eastward (E).  (A point of the form $(i',\infty)$ can only be reached
by ending $P$ with  an infinite number of northward steps along the
line $x=i'$.)  If $P$ starts at $(i,j)$, then we label an eastward
step along the line $y=j'$ with the label
$$
L(s_t)=(j'-j+1)^{(k)}
$$
for some $k$ which can vary with the step, $1\le k\le n$.  
Considering  $P$ as a multiset of labels, it has a well-defined
multidegree $\vec{t}$.
Then for $\vec{t}\vdash t$ we have
$$
h_{\vec{t}}=\sum_P M_P\qmq{where}
M_P=\prod_{\ell^{(k)}}x_\ell^{(k)},    
$$
the sum being over all paths of multidegree $\vec{t}$ from $(i,j)$ to
$(i+t,\infty)$ and the 
product being over all labels in $P$ taken with multiplicity.
Note also that if the labels in 
$P$ are read off from 
left to right, then they correspond to a single-rowed dotted Young
tableau of multidegree $\vec{t}$.

To get products of complete homogeneous symmetric functions and
tableaux of shape $\la=(\la_1,\ldots,\la_l)$,  consider initial vertices
$u_1,\ldots,u_l$  and final vertices $v_1,\ldots,v_l$ with coordinates
\beq
\label{uv}
u_i=(-i,1)\qmq{and} v_i=(\la_i-i,\infty)
\eeq
for $1\le i\le l$.  Consider a family of labeled paths
$\cP=(P_1,\ldots,P_l)$ where, for each $i$, $P_i$ is a path from $u_i$
to $v_{g(i)}$ for some $g\in\fS_l$.  We assign to $\cP$ a {\it
monomial\/} and a {\it sign\/} by
$$
M_\cP=\prod_{i=1}^l M_P\qmq{and} (-1)^\cP=(-1)^g,
$$
respectively.  So denoting the determinant by $D$, we have 
$$
D=\sum_{\cP} (-1)^\cP M_\cP
$$
where the sum is over all path families with beginning and ending
vertices given by~\ree{uv}.

Construct a monomial-preserving, sign-reversing involution $\imath$ on such
$\cP$ which are intersecting as follows. 
Let $i$ be the smallest index such that
$P_i$ intersects some $P_j$ and take $j$ minimum.  Consider the NE-most
point, $v_0$, of $P_i\cap P_j$.   Create $\cP'=\imath \cP$ by replacing
$P_i,P_j$ with $P_i',P_j'$, respectively, where $P_i'$ goes
from $u_i$ to $v_0$ along $P_i$ and then continues along $P_j$, and
similarly for $P_j'$.

Because $\imath$ pairs up intersecting path families of the same
monomial and opposite sign, they all cancel from the determinant
leaving only nonintersecting families.  Furthermore, by the choice of
initial and final points, a family can only be nonintersecting if its
associated element of $\fS_l$ is the identity.  So we now have
$$
D=\sum_{\cP} M_\cP
$$
where the sum is over all nonintersecting families.  But there is a
bijection between such families and tableaux.  Given $\cP$, read the
elements of $P_i$ from left to right to obtain the $i$th row of the
associated tableau $\Td$.  The fact that $\cP$ is nonintersecting is
equivalent to the fact that $\Td$ has increasing columns.  The given
initial and final vertices ensure that the shape of $\Td$ is $\la$.
Applying $\spn{\vec{m}}$ to both sides of the last equality restricts
the multidegree so as to finish the proof of the theorem.\qed

Note that if $\vec{t}=[t]$ has a single component, then 
$\vec{t}\ptn \la_i-i+j$ forces $t=\la_i-i+j$.  So the sums in  
each entry of the determinants reduce to a single term and we recover
the ordinary form of Jacobi-Trudi.

We now specialize to the case of noncommuting variables case so as to
determine the image of $S_\la$ under the involution $\omega$.
\bco
We have $\omega(S_\la)=S_{\la'}$.
\eco
\proof
Merely note that $\omega$ exchanges the two Jacobi-Trudi
determinants.\qed


\section{The Robinson-Schensted-Knuth map}

In this section, we give a generalization of the famous
Robinson-Schensted-Knuth bijection~\cite{knu:pmg,rob:rsg,sch:lid} to
tableaux of arbitrary multidegree.  

A {\it biword of length $n$ over $A$\/} is a $2\times n$ array $\be$ of
elements of $A$ such that if the dots are removed then the columns are
ordered lexicographically with the top row taking precedence. The lower
and upper rows of $\be$ are denoted $\bec$ and $\beh$, respectively.
Viewing $\bec$ and $\beh$ as multisets,  the multidegree of $\be$ is
the pair
$$
\vm(\be)=(\vm(\bec),\vm(\beh)).
$$

We now define a map $\be\stackrel{\rm R-S-K}{\mapsto}(\Td,\Ud)$ whose
image is all pairs of dotted semistandard Young Tableaux of the same shape.
Peform the ordinary Robinson-Schensted-Knuth algorithm on $\be$
(see~\cite{sag:sg} for an exposition) by merely ignoring the dots and
just having them ``come along for the ride.''  For example,  if
$$
\be=\barr{cccccc}
\oned	&\twod	&\twodd	&\twod	&\thrdd	&\foud\\
\twod	&\onedd	&\thrdd	&\thrd	&\twodd	&\oned
\earr
$$
then the sequence of tableaux built by the algorithm is
$$
\begin{array}{lllllll}
\twod\  ,&\onedd\ ,&\onedd\ \thrdd\ ,&\onedd\ \thrdd\ \thrd\ ,
		&\onedd\ \twodd\ \thrd\ ,&\onedd\ \oned\  \thrd\  &\\
	 &\twod	   &\twod	     &\twod	
		&\twod\	\thrdd		 &\twod\  \twodd	  & = \Td,\\
	 &	   &		     &
		&			 &\thrdd		  &\\
&&&&&&\\
\oned\  ,&\oned\  ,&\oned\  \twodd\ ,&\oned\  \twodd\ \twod\ ,
		&\oned\  \twodd\ \twod\ ,&\oned\  \twodd\ \twod\  &\\
	 &\twod	   &\twod	     &\twod	
		&\twod\	\thrdd		 &\twod\  \thrdd	  & = \Ud.\\
	 &	   &		     &
		&			 &\foudd		  &\\

\earr
$$
The next theorem follows directly from the definitions and the
analogous result for the ordinary Robinson-Schensted-Knuth map.
\bth
\label{rsk}
The map 
$$
\be\stackrel{\rm R-S-K}{\mapsto}(\Td,\Ud)
$$
is a bijection between biwords and pairs of dotted semistandard Young
tableaux of the same shape such that
$$
\vm(\be)=(\vm(\Td),\vm(\Ud)).\qed
$$
\eth

Because this analogue is so like the original, most of the properties
of the ordinary Robinson-Schensted-Knuth correspondence carry over
into this setting with virtually no change.  By way of illustration,
here is the corresponding Cauchy identity~\cite{lit:tgc} which follows
directly by 
turning each side of the previous bijection into a generating
function.  Note that for $\beh$ and $\Ud$ we are using a second set of
variables $\byd,\bydd,\ldots,\byddd$.
\bth
We have
$$
\sum_{m\ge0}\sum_{\la,\vm,\vec{p}\ \vdash m}
S_\la^{\vm}(\bxd,\ldots,\bxddd)S_\la^{\vec{p}}(\byd,\ldots,\byddd)=
\prod_{i,j\ge1}
\frac{1}{1-\sum_{k,l=1}^n x_i^{(k)} y_j^{(l)}}.\qed
$$
\eth

\section{Comments and questions}

\quad (I)  Rosas~\cite{ros:sms} computed specializations of symmetric
functions in noncommuting variables and, more generally, of MacMahon symmetric
functions.  

\medskip

\quad (II) Is there an expression for $S_\la^{\vec{m}}$ analogous to
Jacobi's bialternant formula~\cite{mac:sfh,sag:sg,sta:ec2}?

\medskip

\quad (III) Is there a connection between $\Pi(\bx)$ and the partition algebra
$P_n(x)$~\cite{dw:par,hal:cpa,hf:cop,mar:spa}?

\medskip

\quad (IV)  Is there a way to define functions $S_\pi$ for {\it set
partitions\/} $\pi\vdash[n]$ having properties analogous to the
ordinary Schur functions $s_\la$?

\medskip

\quad (V) Given a basis $b_\la(\bx)$ for $\La(\bx)$ we say that
$f(\bx)\in\La(\bx)$ is {\it $b$-positive\/} if the coefficients in the
expansion $f(\bx)=\sum_\la c_\la b_\la(\bx)$ satisfy $c_\la\ge0$ for all $\la$.
Stanley~\cite{sta:sfg,sta:gcr} showed that associated with any
combinatorial graph $G$ there is a symmetric function
$X_G(\bx)\in\La(\bx)$ which generalizes the chromatic polynomial of $G$.
Together with Stembridge~\cite{ss:ijt}, he conjectured that for a
certain family $\cG$ of graphs (those associated with 
({\bf 3+1})-free posets) $X_G$ is $e$-positive for all $G\in\cG$.
Gasharov~\cite{gas:igt} has proved the weaker 
result that  $X_G$ is $s$-positive for all $G\in\cG$.  Gebhard
and Sagan~\cite{gs:csf} have proved that $X_G$ is $e$-positive for all
$G$ in a subfamily of $\cG$ by using symmetric functions in
noncommuting variables.  It would be interesting to enlarge the
subfamily to which these methods can be applied. 

\bigskip

\noindent{\it Acknowledgments.}  Part of the research for this paper
was done while both authors were resident at the Isaac Newton
Institute for Mathematical Sciences in Cambridge, England.  We would
like to thank the Institute 
for support during this period.  We are also indebted to Timothy Chow,
Ira Gessel, and Larry Smith for helpful discussions.  Finally, we owe
a great debt of thanks to an anonymous referee whose incredibly detailed and
helpful comments greatly improved the exposition.

\begin{\bib}{99}

\bibitem{ani:hit} D. J. Anick, On the homogeneous invariants of a
tensor algebra, in ``Algebraic Topology: Proceedings of the
International Conference held March 21--24, 1988,'' Mark Mahowald and
Stewart Priddy eds., Contemporary Mathematics, Vol.\ 96, American
Math.\ Society, Providence, RI, 1989, 15--17.

\bibitem{bc:sef} G. M. Bergman and P. M. Cohn, Symmetric elements in
free powers of rings, {\it J. London Math.\ Soc.\ (2)} {\bf 1} (1969)
525--534.   

\bibitem{dw:par} W. Doran and D. Wales, The partition algebra
revisited, {\it J. Algebra\/} {\bf 231} (2000), 265--330.

\bibitem{dou:sft} P. Doubilet, On the foundations of combinatorial
theory.  VII: Symmetric functions through the theory of distribution
and occupancy, {\it Studies in Applied Math.\/} {\bf 51} (1972),
377--396. 

\bibitem{fg:nsf} S. Fomin and C. Greene, Noncommutative Schur
functions and their applications, {\it Discrete Math.\/}
{\bf 193}  (1998), 179--200.

\bibitem{ful:yt} W. Fulton, ``Young Tableaux,'' London Mathematical
Society Student Texts 35, Cambridge University Press, Cambridge, 1999.

\bibitem{gas:igt} V. Gasharov, Incomparability graphs of ({\bf 3}+{\bf
1})-free posets are  $s$-positive, {\it Discrete Math.} {\bf 157}
(1996), 193--197.

\bibitem{gs:csf} D. Gebhard and B. Sagan, A chromatic symmetric
function in noncommuting variables, {\it J. Algebraic Combin.} {\bf
13} (2001), 227--255. 

\bibitem{gkllrt:nsf} I. M. Gelfand, D. Krob, A. Lascoux, B. Leclerc,
V. Retakh, J.-I. Thibon,  Noncommutative symmetric functions, {\it
Adv.\ in Math.\/} {\bf 112} (1995) 218--348. 

\bibitem{gv:bdp} I.  Gessel  and  G.  Viennot,  Binomial 
determinants, paths, and hook length formulae, {\it \aim} {\bf 58}
(1985), 300--321.

\bibitem{hal:cpa} T. Halverson, Characters of the partition algebra,
{\it J. Algebra\/} {\bf 238} (2001), 502--533.

\bibitem{hf:cop} T. Halverson and J. Fraina, Character orthogonality
fro the partition algebra and fixed points of permutations.
{\it Adv.\ Appl.\ Math.\/} {\bf 31} (2003), 113--131. 

\bibitem{kha:aif} V. K. Kharchenko, Algebras of invariants of free
algebras, {\it Algebra i Logika\/} {\bf 17} (1978) 478--487 (Russian);
{\it Algebra and Logic\/} {\bf 17} (1978), 316--321 (English translation).

\bibitem{knu:pmg} D. E. Knuth, Permutations, matrices and  generalized  Young 
tableaux, {\it Pacific J. Math.} {\bf 34} (1970), 709--727.

\bibitem{ls:mp} A. Lascoux and M.-P. Sch\"utzenberger,
Le monoid plaxique, in ``Noncommutative Structures in Algebra and
Geometric Combinatorics, (Naples, 1978),''  Quad.\ Ricerca Sci.,
Vol.\ 109, CNR, Rome, 1981, 129--156.

\bibitem{lin:vri} B. Lindstr\"om, On the vector representation of induced
matroids, {\it Bull.\  London Math.\  Soc.\ } {\bf 5} (1973), 85--90.

\bibitem{lit:tgc} D. E. Littlewood, ``The Theory of Group Characters,''
Oxford University Press, Oxford, 1950.

\bibitem{mac:sfh} I. G. Macdonald,  ``Symmetric  functions 
and Hall polynomials,'' 2nd edition, \oup, Oxford, 1995.

\bibitem{mac:ca} P. A. MacMahon, ``Combinatorial Analysis,''
Vols. 1 and 2, Cambridge University Press, Cambridge, 1915, 1916; reprinted by
Chelsea, New York, NY, 1960.

\bibitem{mar:spa} P. Martin, The structure of partition algebras, 
{\it J. Algebra\/} {\bf 183} (1996) 319--358.

\bibitem{nw:wgp} S. D. Noble and D. J. A. Welsh, A weighted graph
polynomial from chromatic invariants of knots, 
Symposium \`a la M\'emoire de Fran\c{c}ois Jaeger (Grenoble, 1998) 
{\it Annales l'Institut Fourier} {\bf 49} (1999),
1057--1087.

\bibitem{rob:rsg} G. de B. Robinson, On  representations  of  the  symmetric 
group, {\it Amer. J. Math.} {\bf 60} (1934), 745--760.

\bibitem{ros:msf} M. H. Rosas, MacMahon symmetric functions, the
partition lattice, and Young subgroups, {\it J. Combin.\ Theory Ser.\ A}
{\bf 96} (2001), 326--340.

\bibitem{ros:sms} M. H. Rosas, Specializations of MacMahon symmetric
functions and the polynomial algebra,  {\it Discrete Math.\ }
{\bf 246} (2002), 285--293.

\bibitem{rot:tmf} G.-C. Rota, On the foundations of combinatorial
theory I. Theory of M\"obius functions, {\it Z. 
Wahrscheinlichkeitstheorie} {\bf 2} (1964), 340--368.

\bibitem{sag:sg} B. Sagan, ``The Symmetric Group: Representations,
Combinatorial Algorithms,  and Symmetric Functions,'' 2nd edition,
Springer-Verlag, New York, 2001. 

\bibitem{sch:lid}  C.   Schensted,   Longest   increasing   and   decreasing 
subsequences, {\it Canad. J. Math.} {\bf 13} (1961), 179--191.

\bibitem{sch:ukm} I. Schur, ``\"{U}ber eine Klasse von Matrizen die 
sich einer gegebenen Matrix zuordnen lassen,'' Inaugural-Dissertation, Berlin,
1901.

\bibitem{sta:sfg} R. P. Stanley, A symmetric function generalization
of the chromatic polynomial of a graph, {\it Advances in Math.}\ 
{\bf 111} (1995), 166--194.

\bibitem{sta:ec1} R. P. Stanley, ``Enumerative Combinatorics,
Volume 1,''  Cambridge University Press, Cambridge, 1997.

\bibitem{sta:gcr} R. P. Stanley, {\bf G}raph Colorings {\bf a}nd
{\bf r}elated {\bf s}ymmetric functions: {\bf i}deas and 
{\bf a}pplications: A description of results, interesting applications,
\& notable open problems,  Selected papers in honor of Adriano Garsia
(Taormina, 1994) {\it Discrete Math.\/} {\bf 193} (1998), 267--286.

\bibitem{sta:ec2} R. P. Stanley, ``Enumerative Combinatorics,
Volume 2,''  Cambridge University Press, Cambridge, 1999.

\bibitem{ss:ijt} R. P. Stanley and J. Stembridge, On immanants of
Jacobi-Trudi  matrices and permutations with restricted position, 
{\it J. Combin.\ Theory Ser.\ A}  {\bf 62} (1993), 261--279.

\bibitem{wei:ati} L. Weisner, Abstract theory of inversion of finite
series, {\it Trans.\  Amer.\  Math.\  Soc.\ } {\bf 38} (1935) 474--484.

\bibitem{whi:lem}  H. Whitney, A logical expansion in mathematics,
{\it Bull.\ Amer.\ Math.\ Soc.\ } {\bf 38} (1932), 572--579.

\bibitem{wol:sfn} M. C. Wolf,  Symmetric functions of noncommuting
elements,  {\it Duke Math.\ J.} {\bf 2} (1936) 626--637.

\end{\bib}

\end{document}